\documentclass[a4paper]{amsart}
\usepackage{a4}

\title{A conjecture of Bavula on homomorphisms of the Weyl algebra}

\newcommand{\de}{\partial}

\newcommand{\ad}{\mbox{ad}}
\newcommand{\Ad}{\mbox{Ad}}

\newcommand{\Z}{\Bbb Z}
\newcommand{\J}{\mbox{J}}

\newcommand{\Def}{\mbox{def}}

\newcommand{\gkdim}{\mbox{GKdim}}

\begin{document}

\date{}
\maketitle

\begin{center}
\textbf{Leonid Makar-Limanov} \footnote{Supported
by an NSA grant H98230-09-1-0008, by an NSF grant DMS-0904713, and a Fulbright fellowship awarded by the United States--Israel Educational Foundation.\\
For the reader's convenience the proofs of all necessary lemmas are included, eliminating the need to look for them elsewhere.}\\
Department of Mathematics \& Computer Science,\\
the Weizmann Institute of Science, Rehovot 76100, Israel and\\
Department of Mathematics, Wayne State University, \\
Detroit, MI 48202, USA\\
e-mail: {\em lml@math.wayne.edu}\\
\end{center}

\begin{abstract}
In the paper {\em The inversion formulae for automorphisms of
polynomial algebras and differential operators in prime
characteristic}, J. Pure Appl. Algebra 212 (2008), no. 10,
2320--2337, see also Arxiv:math.RA/0604477, Vladimir Bavula states
the following Conjecture:

(BC) Any endomorphism of a Weyl algebra (in a finite characteristic
case) is a monomorphism.

The purpose of this preprint is to prove BC for $A_1$, show that BC
is wrong for $A_n$ when $n > 1$, and prove an analogue of $BC$ for symplectic
Poisson algebras.
\end{abstract}

The Weyl algebra $A_n$ is an algebra over a field $F$ generated by
$2n$ elements $x_1, \ldots x_n; \ y_1, \ldots, y_n$ which satisfy
relations $[x_i, y_j] (= x_iy_j - y_jx_i) = \delta_{ij}, \ \
[x_i,x_j]=0, \ \ [y_i,y_j] = 0, $ where $\delta_{ij}$ is the
Kronecker symbol and $1\leq i,j\leq n$. Weyl algebras appeared quite
some time ago and initially were considered only over fields of
characteristic zero. Arguably the most famous algebraic problem
related to these algebras is the Dixmier conjecture (see [D]) that
any homomorphism of $A_n$ is an automorphism if char$(F) = 0$. This
problem is still open even for $n = 1$.

The finite characteristic case is certainly less popular but lately
appears to attract more attention because it helps to connect
questions related to the Weyl algebras and to polynomial rings,
e. g. to connect the famous Jacobian Conjecture with the Dixmier
conjecture (see [T1], [BK], and [AE]). There is a striking difference
between the zero characteristic and the finite characteristic cases.
While for characteristic zero the center of $A_n$ is just the ground
field and $A_n$ is infinite-dimensional over the center, when the
characteristic is not zero the center of $A_n$ is a polynomial ring
in $2n$ generators and $A_n$ is a finite-dimensional free module over the center.\\

A vector space $B$ with two bilinear operations
$x\cdot y$ (a multiplication) and $\{x,y\}$ (a Poisson bracket) is
called {\em a Poisson algebra} if $B$ is a commutative associative
algebra under $x\cdot y$, $B$ is a Lie algebra under $\{x,y\}$, and
$B$ satisfies the Leibniz identity: $ \{x, y\cdot z\}=\{x,y\}\cdot z
+ y\cdot \{x,z\}. $ \\

Here we will be concerned with symplectic (Poisson) algebras $PS_n$.
For
each $n$ the algebra $PS_n$ is a polynomial algebra $F[x_1, \ldots x_n; \
y_1, \ldots, y_n]$ with the Poisson bracket defined by $
\{x_i,y_j\}=\delta_{ij}, \ \ \{x_i,x_j\}=0, \ \ \{y_i,y_j\}=0, $
where $\delta_{ij}$ is the Kronecker symbol and $1\leq i,j\leq n$.
Hence $\{f, g\} = \sum_i({\de f\over \de x_i}{\de g\over \de y_i} -
{\de f\over \de y_i}{\de g\over \de x_i})$.

To distinguish $PS_n$ and $A_n$ we will write $PS_n$ as $F\{x_1,y_1,
\ldots,x_n,y_n\}.$ One can think about $PS_n$ as a commutative
approximation of $A_n$ (and of $A_n$ as a quantization of $PS_n$).\\

It is clear that $PS_n$ is a polynomial algebra with some additional
structure. It is less clear what is $A_n$. Of course, we can think
about a Weyl algebra as a factor algebra of a free associative
algebra by the ideal $I$ which  corresponds to the relations, but it
is not obvious even that $1 \not\in I$ so the resulting factor
algebra may be the zero algebra.\\

\textbf{Lemma on basis.} The monomials $y_1^{j_1}x_1^{i_1} \ldots
y_n^{j_n}x_n^{i_n}$ form a basis $\mathcal{P}$ of $A_n$ over $F$.\\
\textbf{Proof.} Any monomial $\mu$ in $A_n$ (which in this
consideration may be the zero algebra) can be written as a product
$\mu = \mu_1 \mu_2 \ldots \mu_n$ where $\mu_i$ is a monomial of
$x_i, \ y_i$ since different pairs $(x_j, y_j)$ commute.
Furthermore, since $x_iy_i = y_ix_i + 1$ any monomial in $x_i, \
y_i$ can be written as a linear combination of monomials
$y_i^kx_i^l$ with coefficients in $\Z$ or in $\Z_p$ depending on the
characteristic of $F$.

It remains to show that the monomials in $\mathcal{P}$ are linearly
independent over $F$. This can be done by finding a homomorphic
image of $A_n$ in which the images of monomials from $\mathcal{P}$
are linearly independent.

If char$F = 0$ there is a natural representation of $A_n$. Consider
homomorphisms $X_i$ and $Y_j$ of $R = F[y_1, \ldots, y_n]$ defined
by $X_i(r) =  \frac{\de \, r}{\de \, y_i}$ and $Y_j(r) = y_jr$. A
straightforward computation shows that  $\alpha(x_i) = X_i$,
$\alpha(y_j) = Y_j$ defines a homomorphism of $A_n$ into the ring of
homomorphisms of $R$ and that the images of monomials from
$\mathcal{P}$ are linearly independent. Unfortunately this
representation is not satisfactory when char$F = p \neq 0$. In that case it is
easy to see that $X_i^p = 0$.

Here is a way to remedy the problem. Consider $R_n = R[z_1, \ldots,
z_n]$, and homomorphisms $X_i, \ Y_j$, and $Z_k$ of $R_n$ defined by
$X_i(r) =  \frac{\de \, r}{\de \, y_i}$, $Y_j(r) = y_jr$ and $Z_k(r)
= z_kr$ for $r \in R_n$. Since $Z_k$ commute with $X_i$ and  $Y_j$,
if we replace $X_i$ by $\widehat{X_i} = X_i + Z_i$ then
$[\widehat{X_i}, Y_j] = [X_i, Y_j]$ and $[\widehat{X_i}, \widehat{X_j}] = 0$. Now, for $\sigma = \sum
f_{\mathbf{i}\mathbf{j}} Y_1^{j_1}\widehat{X}_1^{i_1} \ldots
Y_n^{j_n}\widehat{X}_n^{i_n}$ we have $\sigma(1) = \sum
f_{\mathbf{i}\mathbf{j}} y_1^{j_1}z_1^{i_1} \ldots
y_n^{j_n}z_n^{i_n}$, which is zero only if all
$f_{\mathbf{i}\mathbf{j}} = 0$. (Here $\mathbf{i}$ and $\mathbf{j}$ are multi-indices as usual.) $\Box$\\

Let us call the presentation of an element $a \in A_n$ through the
basis $\mathcal{P}$ \textsl{standard}. Further we will use only the
standard presentations of elements of $A_n$. So $A_n$
is isomorphic to a corresponding polynomial ring as a vector space.\\

\textbf{Remark 1.} $A_n$ is a domain (does not have zero divisors).
To see this consider a degree-lexicographic ordering of monomials in
$\mathcal{P}$ first by the total degree and then by $y_1
>> x_1 >> y_2 >> x_2 \ldots
>> y_n >> x_n$. Then the commutation relations of $A_n$ give
$|fg| = |f||g||$ where $|h|$ for $h \in A_n \setminus 0$ denotes the
largest monomial appearing in $h$.  $\Box$ \\

If char$(F) = 0$ then BC is very easy to prove (and is well-known)
both for $A_n$ and $PS_n$. Suppose that $\varphi$ has a non-zero
kernel. Let us take a non-zero element in the kernel of $\varphi$ of
the smallest total degree $\deg$ possible. It is clear that $\deg(ab) = \deg(a) + \deg(b)$ for
$a, \ b \in A_n$ because of the commutation relations.
Consider $PS_n$ first. If $\Lambda$ is a ``minimal" element of the
kernel then $\{x_i, \Lambda\} = {\de \Lambda\over \de y_i}$ and
$\{y_i, \Lambda\} = -{\de \Lambda\over \de x_i}$ should be
identically zero because of the minimality of $\Lambda$. So ${\de
\Lambda\over \de y_i} = 0$ and ${\de \Lambda\over \de x_i} = 0$ for
all $i$. If char$(F) = 0$ this means that $\Lambda \in F$. But our
homomorphism is over $F$, so $\Lambda = 0$. A similar proof works for
$A_n$ where the elements $[x_i, \Lambda]$ and $[y_i, \Lambda]$
should be identically zero which again shows that $\Lambda = 0$.\\

\textit{From now on we assume that char$(F) = p \neq 0$.} \\

Let us start with purely computational observations.

A straightforward computation shows that $[ab, c] = [a, c]b  + a[b,
c]$. Therefore $[x_1^{k+1}, y_1] = [x_1^k, y_1]x_1 + x_1^k[x_1,
y_1]$ and since $[x_1, y_1] = 1$ induction on $k$ gives $[x_1^k,
y_1] = kx_1^{k-1}$. Similarly, $[x_1, y_1^k] = ky_1^{k-1}$ and the index
1 can be replaced by any $i \in \{1, 2, \ldots, n\}$.

Denote $[a, b]$ by $\ad_a(b)$. We will use the following equality: $\ad_a^p(B)
= \ad_{a^p}(b)$. In order to prove it observe that $\ad_a(b) = (a_l
- a_r)(b)$ where $a_l$ and $a_r$ are the operators of left and right
multiplication by $a$. It is clear that $a_l$ and $a_r$ commute. So
$\ad_a^p(b) = (a_l - a_r)^p(b) = (a_l^p -
a_r^p)(b) = \ad_{a^p}(b)$.\\

\textbf{Lemma on center.} (a) The center $Z(A_n)$ of $A_n = F[x_1,
\ldots x_n; \ y_1, \ldots, y_n]$ is the polynomial ring $F[x_1^p,
\ldots x_n^p; \ y_1^p, \ldots, y_n^p]$.\\
(b) The Poisson center of $PS_n = F\{x_1, \ldots x_n; \ y_1, \ldots,
y_n\}$ is the polynomial ring $F[x_1^p, \ldots x_n^p; \ y_1^p,
\ldots, y_n^p]$.\\
\textbf{Proof.} (a) Consider $x_1^p$. It is clear from the
definition of $A_n$ that $x^p_1$ commutes with all generators with a
possible exception of $y_1$. As we observed above, $[x_1, y_1] = 1$
implies $[x_1^k, y_1] = kx_1^{k-1}$. So $[x_1^p, y_1] = px_1^{p-1} =
0$ and $x^p_1$ is in the center of $A_n$. Similarly all $x_j^p$ and
$y_j^p$ are in the center and $Z(A_n) \supseteq E$ where $E =
F[x^p_1, \ldots x^p_n; \ y^p_1, \ldots, y^p_n]$.

Any element $a \in A_n$ can be written as $a = \sum c_{\mathbf{i},
\mathbf{j}} y_1^{j_1} x_1^{i_1} \ldots y_n^{j_n} x_n^{i_n}  = a_0 +
\sigma$ where $0 \leq i_s < p$ and $0 \leq j_s < p$, $c_{\mathbf{i},
\mathbf{j}} \in E$, $a_0 \in E$, and $\sigma$ is the sum of all
monomials of $a$ which do not belong to $E$. If $a \in Z(A_n)$ then
$[x_1, a] = 0$. Now, $[x_1, y_1^{j_1} x_1^{i_1}  \ldots y_n^{j_n}
x_n^{i_n}] = [x_1, y_1^{j_1}] x_1^{i_1} \ldots y_n^{j_n} x_n^{i_n} =
j_1y_1^{j_1-1}x_1^{i_1} \ldots y_n^{j_n} x_n^{i_n}$ and we have
similar formulae when $x_1$ is replaced by any $x_i$ or $y_j$. So if
one of the monomials in $\sigma$ is not zero we can take the commutator of $a$ with an
appropriate $x_i$ or $y_j$ and obtain a non-trivial linear
dependence between monomials of $\mathcal{P}$ contrary to the Lemma on
basis.\\
(b) An element $a$ belongs to the Poisson center
$Z(\mathcal{A})$ of a Poisson algebra $\mathcal{A}$ if $\{a, b\} =
0$ for all $b \in \mathcal{A}$. If $f \in Z(PS_n)$ then ${\de f\over
\de x_i} = \{f, y_i\} = {\de f\over \de y_i} = \{x_i, f\} = 0$ for
all $i$ which is possible only if $f \in F[x_1^p, \ldots x_n^p; \
y_1^p, \ldots,
y_n^p]$.  $\Box$\\

\textbf{Nousiainen Lemma.} Let $\varphi$ be a homomorphism of $A_n$
or $PS_n$ correspondingly. Then (a) $A_n = Z(A_n)[\varphi(x_1),
\ldots \varphi(x_n); \ \varphi(y_1), \ldots, \varphi(y_n)]$;\\ (b)
$PS_n = Z(PS_n)[\varphi(x_1), \ldots \varphi(x_n); \ \varphi(y_1),
\ldots, \varphi(y_n)]$.\\
\textbf{Proof.} (a) Let $E = F[x^p_1, \ldots x^p_n; \ y^p_1, \ldots,
y^p_n]$. From the Lemma on center $Z(A_n) = E$. The algebra $A_n$ is a
finite-dimensional module over $E$ since any element $a \in A_n$ can
be written as $a = \sum c_{\mathbf{i}, \mathbf{j}} y_1^{j_1}
x_1^{i_1}  \ldots y_n^{j_n} x_n^{i_n} $ where $0 \leq i_s < p$, $0
\leq j_s < p$, and  $c_{\mathbf{i}, \mathbf{j}} \in E$. Let $K =
F(x^p_1, \ldots x^p_n; \ y^p_1, \ldots, y^p_n)$ be the field of
fractions of $E$ and let $D_n = K[x_1, \ldots x_n; \ y_1, \ldots,
y_n]$. Algebra $D_n$ is a skew-field since $D_n$ is a
finite-dimensional vector space over $K$ and $D_n$ does not have
zero divisors according to Remark 1. (Recall that $x_i^p, \ y_j^p
\in K$, so any $f \in D_n$ satisfies a non-zero relation
$\sum_{i=0}^N k_if^i = 0$ where $k_i \in K$ and $N \leq p^{2n}$.)

The monomials $y_1^{j_1} x_1^{i_1}  \ldots y_n^{j_n} x_n^{i_n} $, where
$0 \leq i_s < p$, $0 \leq j_s < p$, are linearly independent over
$K$. Indeed, if $\Lambda = \sum c_{\mathbf{i}, \mathbf{j}} y_1^{j_1}
x_1^{i_1}  \ldots y_n^{j_n} x_n^{i_n}  = 0$ where
$c_{\mathbf{i}, \mathbf{j}} \in K$ then $[x_m,  \Lambda] =
[y_m, \Lambda] = 0$ and we can obtain from a non-trivial dependence
$\Lambda$ a ``smaller" one. So we can invoke induction on, e. g. the
sum of the total degrees of monomials in $\Lambda$.

Since any element $a \in D_n$ can be written as $a = \sum
c_{\mathbf{i}, \mathbf{j}} y_1^{j_1} x_1^{i_1}  \ldots y_n^{j_n}
x_n^{i_n} $ where $0 \leq i_s < p$, $0 \leq j_s < p$, and
$c_{\mathbf{i}, \mathbf{j}} \in K$, the dimension of $D_n$ over $K$
is $p^{2n}$ .

Consider now monomials $v_1^{j_1} u_1^{i_1} \ldots v_n^{j_n}
u_n^{i_n}$ where $u_i = \varphi(x_i)$, $v_j =  \varphi(y_j)$, $0
\leq i_s < p$, and $0 \leq j_s < p$. Let us check that they are also
linearly independent over $K$. If $\Lambda = \sum c_{\mathbf{i},
\mathbf{j}} v_1^{j_1} u_1^{i_1} \ldots v_n^{j_n} u_n^{i_n} = 0$ then
$[u_m,  \Lambda] = [v_m, \Lambda] = 0$ and, since the commutation
relations are the same as above, we obtain from a non-trivial
dependence $\Lambda$ a ``smaller" one.

Since there are exactly $p^{2n}$ of these monomials,
they also form a basis of $D_n$ over $K$
and
any element $a \in A_n \subset D_n$ can be written as $a = \sum
c_{\mathbf{i}, \mathbf{j}} v_1^{j_1} u_1^{i_1} \ldots v_n^{j_n}
u_n^{i_n}$ where
$c_{\mathbf{i}, \mathbf{j}} \in K$.

It remains to show that all
$c_{\mathbf{i}, \mathbf{j}} \in E$.
Order the monomials $v_1^{j_1} u_1^{i_1} \ldots v_n^{j_n} u_n^{i_n}$
degree-lexicographically as in Remark 1. Let $\mu = v_1^{j_1}
u_1^{i_1} \ldots v_n^{j_n} u_n^{i_n}$ be the largest monomial. Then
$\ad_{v_n}^{i_n}\ad_{u_n}^{j_n}  \ldots
\ad_{v_1}^{i_1}\ad_{u_1}^{j_1}(a) = (-1)^I \prod_{m=1}^n (i_m)!
(j_m)! c_{\mathbf{i}, \mathbf{j}}$, where $I = \sum_{m=1}^n
i_m$, belongs to $A_n$. Since $(-1)^I \prod_{m=1}^n (i_m)! (j_m)! \neq 0$ we conclude
that $c_{\mathbf{i}, \mathbf{j}} \in A_n \bigcap K = E$, replace $a$
by $a - c_{\mathbf{i}, \mathbf{j}} \mu$, and finish
by induction on the number of monomials of $a$. \\
(b) Let $T_n$ be
the field of rational functions $F(x_1, \ldots x_n; \ y_1, \ldots,
y_n )$ endowed with the same bracket as $PS_n$: $\{f, g\} =
\sum_i({\de f\over \de x_i}{\de g\over \de y_i} - {\de f\over \de
y_i}{\de g\over \de x_i})$. Then $T_n$ becomes a Poisson algebra and
we can think of $PS_n$ as a subalgebra of $T_n$. It is clear that
$Z(T_n) = F(x_1^p, \ldots x_n^p; \ y_1^p, \ldots, y_n^p)$ and that
$T_n$ is a $p^{2n}$-dimensional vector space over $Z(T_n)$.

Denote $u_i = \varphi(x_i), \ v_i = \varphi(y_i)$. To show that the
monomials $v_1^{j_1} u_1^{i_1} \ldots v_n^{j_n} u_n^{i_n}$ where $0
\leq i_s < p$ and $0 \leq j_s < p$ are linearly independent over
$Z(T_n)$ we, as above, can consider a ``minimal" relation $\Lambda =
\sum c_{\mathbf{i}, \mathbf{j}} v_1^{j_1} u_1^{i_1} \ldots v_n^{j_n}
u_n^{i_n} = 0$ and get ``smaller" relations by taking $\{u_m,
\Lambda\}$ and $\{v_m, \Lambda\}$.

If $a \in PS_n$ is presented as $a
=  \sum c_{\mathbf{i}, \mathbf{j}} v_1^{j_1} u_1^{i_1} \ldots
v_n^{j_n} u_n^{i_n}$ where $0 \leq i_s < p$, $0 \leq j_s < p$, and
$c_{\mathbf{i}, \mathbf{j}} \in Z(T_n)$ then all $c_{\mathbf{i},
\mathbf{j}} \in Z(PS_n)$; to see this just replace, in the considerations above,
$\ad_z$ by $\Ad_z$ defined by $\Ad_z(b) = \{z, b\}$.       $\Box$ \\

\textbf{Corollary.} There are no homomorphisms from $A_n$ into $A_{n-1}$.\\
\textbf{Proof.} Assume that we have a homomorphism $\phi: A_n
\rightarrow A_{n-1}$. Consider images $u_i = \phi(x_i)$ and $v_i =
\phi(y_i)$. According to Nousiainen Lemma $A_{n-1}$ is a vector
space over the center of $A_{n-1}$ with a basis consisting of
monomials $v_1^{j_1} u_1^{i_1} \ldots v_{n-1}^{j_{n-1}}
u_{n-1}^{i_{n-1}}$, $0 \leq i_s < p, \ 0 \leq j_s < p$.
Therefore $u_n$ and $v_n$ are in the center of $A_{n-1}$ and
hence commute with each other.
$\Box$ \\

\textbf{Remark 2}. This Lemma is similar to a lemma from Pekka
Nousiainen's PhD thesis (Pennsylvania State University, 1981) which
was never published (see [BCW]). Nousiainen proved his Lemma in a
commutative setting for a Jacobian set of polynomials, i. e. he
proved that if $z_1, \ldots, z_n \in F[y_1, \ldots, y_n]$ and the
Jacobian $\J(z_1, \ldots, z_n) = 1$ then $F[y_1, \ldots, y_n] =
F[y_1^p, \ldots, y_n^p; z_1, \ldots, z_n]$. This shows immediately
that $F[y_1, \ldots, y_n] = F[y_1^P, \ldots, y_n^P; z_1, \ldots,
z_n]$ where $P = p^m$ and $m$ is any natural number. Indeed, if say
$y_1 = \sum c_{\mathbf{i}} z_1^{i_1} \ldots z_n^{i_n}$ where
$c_{\mathbf{i}}  \in F[y_1^p, \ldots, y_n^p]$ then $y_1^p = \sum
c_{\mathbf{i}}^p (z_1^{i_1} \ldots z_n^{i_n})^p$ and
$c_{\mathbf{i}}^p \in F[y_1^{p^2}, \ldots, y_n^{p^2}]$.

The same is true for the symplectic algebras. But it is wrong even for $A_1$.

Take e. g. $u = x, \ v = y^2x - y$ when $p = 2$. Then $u^2 = x^2, \
v^2 = y^4x^2$ and $D_1 \neq F(x^4, y^4)[u, v]$ since $u^2$ and $v^2$
are linearly dependent over $F(x^4, y^4)$. Here, of course, $A_1
\neq F[x^4, y^4; u, v]$. This difference between $A_n$ and $PS_n$ is the reason for BC to be correct for $PS_n$ and wrong for $A_n$. 

The Nousiainen Lemma for Weyl algebras is proved in [T2]
and [AE]. $\Box$ \\

\textbf{Theorem 1.} BC is true for Poisson algebras $PS_n$. \\
\textbf{Proof.} In the Nousiainen Lemma we proved that $PS_n =
Z(PS_n)[u_1, \ldots u_n; \ v_1, \ldots, v_n]$ where $u_i =
\varphi(x_i), \ v_i = \varphi(y_i)$. So $a =  \sum c_{\mathbf{i},
\mathbf{j}} v_1^{j_1} u_1^{i_1} \ldots v_n^{j_n} u_n^{i_n}$ where
$c_{\mathbf{i}, \mathbf{j}} \in Z(PS_n) = F[x_1^p, \ldots x_n^p; \
y_1^p, \ldots, y_n^p]$ for any $a \in PS_n$.
Therefore
$$a^p = \sum c_{\mathbf{i}, \mathbf{j}}^p v_1^{pj_1}
u_1^{pi_1} \ldots v_n^{pj_n} u_n^{pi_n} \ {\rm where \ } c_{\mathbf{i},
\mathbf{j}}^p \in F[x_1^{p^2}, \ldots, x_n^{p^2}; \ y_1^{p^2},
\ldots, y_n^{p^2}].$$
Hence \\ $F[x_1^p, \dots, x_n^p; \ y_1^p, \dots,
y_n^p] \subset F[x_1^{p^2}, \dots, x_n^{p^2}; \ y_1^{p^2}, \dots,
y_n^{p^2}][u_1^p, \dots, u_n^p; \ v_1^p, \dots, v_n^p]$
and
$a =
\sum d_{\mathbf{i}, \mathbf{j}} v_1^{j_1} u_1^{i_1} \dots v_n^{j_n}
u_n^{i_n}$ where \\ $d_{\mathbf{i}, \mathbf{j}} \in F[x_1^{p^2}, \dots
x_n^{p^2}; \ y_1^{p^2}, \dots, y_n^{p^2}][u_1^p, \dots u_n^p; \
v_1^p, \ldots, v_n^p].$

So $PS_n = F[x_1^{p^2}, \ldots, x_n^{p^2}; \
y_1^{p^2}, \ldots, y_n^{p^2}][u_1, \ldots, u_n; \ v_1, \ldots, v_n]$.
By iterating this process we will get that $PS_n = F[x_1^P, \ldots
x_n^P; \ y_1^P, \ldots, y_n^P][u_1, \ldots, u_n; \ v_1, \ldots, v_n]$
where $P = p^m$ for any positive integer $m$.

It is clear that $u_i^P, \ v_j^P \in F[x_1^P, \ldots, x_n^P; \ y_1^P,
\ldots, y_n^P]$ so $PS_n$ is spanned over $F[x_1^P, \ldots, x_n^P; \
y_1^P, \ldots, y_n^P]$ by monomials $v_1^{j_1} u_1^{i_1} \ldots,
v_n^{j_n} u_n^{i_n}$ where $0 \leq i_s < P, \ 0 \leq j_s < P$. Of
course, $PS_n$ is spanned over $F[x_1^P, \ldots x_n^P; \ y_1^P,
\ldots, y_n^P]$ by monomials $y_1^{j_1} x_1^{i_1} \ldots y_n^{j_n}
x_n^{i_n}$ where $0 \leq i_s < P, \ 0 \leq j_s < P$ and these
monomials are linearly independent over $F[x_1^P, \ldots x_n^P; \
y_1^P, \ldots, y_n^P]$.

So $F[x_1, \ldots x_n; \
y_1, \ldots, y_n]$ is a free module over $F[x_1^P, \ldots x_n^P; \
y_1^P, \ldots, y_n^P]$ of dimension $P^{2n}$.

If $\varphi$ is not an injection then there is a linear dependence
over $F$ between monomials $v_1^{j_1} u_1^{i_1} \ldots v_n^{j_n}
u_n^{i_n}$ where $0 \leq i_s < P, \ 0 \leq j_s < P$ provided $P$ is
sufficiently large. But this is impossible since these monomials are
linearly independent over $F[x_1^P, \ldots x_n^P; \
y_1^P, \ldots, y_n^P]$.   $\Box$ \\

We prove now using Gelfand-Kirillov dimension that a homomorphism of
$A_1$ into $A_n$ is
an embedding. Here is a definition
of the Gelfand-Kirillov dimension ($\gkdim$) suitable for our
purpose. Let $R$ be a finitely-generated associative algebra over
$F$: $R = F[r_1, \ldots, r_m]$. Consider a free associative algebra
$F_m = F\langle z_1, \ldots, z_m\rangle$ and linear subspaces $F_{m,
N}$ of all elements of $F_m$ with total degree at most $N$. Let
$\alpha$ be a homomorphism of $F_m$ onto $R$ defined by $\alpha(z_i)
= r_i$ and let $R_N = \alpha(F_{m, N})$. Each $R_N$ is a
finite-dimensional vector space (over $F$); denote $d_N =
\dim(R_N)$. $$\gkdim(R) = lim_{N \rightarrow \infty}
\frac{\ln(d_N)}{\ln(N)}.$$

Though this definition uses a particular system of generators, it is
possible to prove that $\gkdim(R)$ does not depend on such a choice
(see [GK] or [KL]). It is not difficult to show that in the
commutative case Gelfand-Kirillov dimension coincides with the
transcendence degree. \\

\textbf{Lemma on GK-dimension.} Let $R = F[z_1, \ldots, z_m]$ be a
ring of polynomials. If $a, b \in R$ are algebraically dependent
then $\gkdim(S) \leq 1$ for any finitely generated subalgebra $S
\subset A = F[a, b]$. \\
\textbf{Proof.} If $a, b \in F$ then $F[a, b] = F$ and any
subalgebra of $A$ is $F$. So in this case $\gkdim(S) = 0$. Assume
now that $a \not\in F$, i. e. $\deg(a) > 0$. If $\deg(b) = 0$ then
$A = F[a]$ and $d_N = N + 1$ where  $d_N = \dim(A_N)$, so $\gkdim(A)
= 1$. Let $S$ be a subalgebra of $A$ generated by $a_1, \ldots,
a_m$. Then $a_i = q_i(a)$ where $q_i$ are polynomials. Assume that
degrees of all these polynomials do not exceed $d$. Therefore a polynomial
$g(a_1, \ldots, a_m)$ of the total degree
$N$ can be rewritten as a polynomial in
$a$ of degree at most $dN$. Hence $N < d_N \leq dN + 1$ if any
of $a_i$ is not in $F$ and $\gkdim(S) = 1$. If all $a_i \in F$ then
$\gkdim(S) = 0$.

Now, let $\deg(b) > 0$. Since $a$ and $b$ are algebraically
dependent, $Q(a, b) = 0$ for a non-zero polynomial $Q$. Order
monomials $a^ib^j$ by total degree $i + j$ and then
lexicographically by $a >> b$. If $\mu = a^kb^l$ is the largest
monomial in $Q$ we can write $\mu = Q_1(a,b)$ where all monomials of
$Q_1$ are less than $\mu$. So we can replace any monomial $\nu =
a^ib^j$  where $i \geq k, \ j \geq l$ by a linear combination of
monomials which are less than $\nu$. Hence any $c \in A$ of the
total degree at most $N$ can be written as a linear combination of
monomials $a^ib^j$ where $i + j \leq N$ and either $i < k$ or $j <
l$. There are less than $(k + l)(N + 1)$ and more than $N$ monomials
satisfying these properties. Therefore  $N < d_N < (k + l)(N + 1)$
and $\gkdim(A) = 1$. If $S$ is a subalgebra of $A$ generated by
$a_1, \ldots, a_m$ then $a_i = q_i(a, b)$ where $q_i$ are
polynomials of total degrees bounded by some $d$.
If we take a polynomial $g(a_1,
\ldots, a_m)$ of total degree $N$ then $g(a_1, \ldots, a_m)$ can
be rewritten as a polynomial in $a$ and $b$ of degree at most $dN$.
Therefore $N < d_N \leq (k + l)(dN + 1)$ if any of $a_i$ is not in $F$ and
$\gkdim(S) = 1$. If all $a_i \in F$ then $\gkdim(S) = 0$. $\Box$\\

Denote by $\deg$ the total degree function on $A_n$ and by
$\overline{a}$ the element of $A$ which is $\deg$ homogeneous and such that
$\deg(a - \overline{a}) < \deg(a)$. We will refer to $\overline{a}$
as the leading form of $a$. From the commutation relations in $A_n$ it follows that
$\overline{ab} = \overline{ba}$ and that $\deg(\overline{[a, b]}) <
\deg(\overline{ab})$.

We will think about leading forms not as elements of $A_n$ but
rather as commutative polynomials. Then $\overline{ab} =
\overline{a} \overline{b} = \overline{b} \overline{a}$.\\

\textbf{Lemma on dependence.} Let $a$ and $b$ be algebraically
dependent non-zero homogeneous polynomials and $q = \deg(a)$, $r
= \deg(b)$. Then $a^r$ and $b^q$ are proportional, i. e. there
exists an $f \in F$ so that $a^r - fb^q = 0$.\\
\textbf{Proof.} The polynomials $a$ and $b$ are algebraically dependent,
so there is a non-zero polynomial $Q$ for which $Q(a, b) = \sum_{i, j}
f_{ij}a^ib^j = 0$. Since $a$ and $b$ are homogeneous we may assume that
$qi + rj$ is the same for all monomials of $Q$. Indeed, any $Q$ can be
presented as $Q = \sum_i Q_k$ where $Q_k$ are $q, r$-homogeneous.
Then either $Q_k(a, b) = 0$ or $\deg(Q_k(a, b)) = k$ and different components cannot cancel out.

Therefore over an algebraic closure $\overline{F}$ of $F$ we can write
$Q = f_0 a^k b^l \prod_i (a^{r'} - f_i b^{q'})$ where $f_i \in \overline{F}$,
$r' = {r\over d}, \ q' = {q \over d}$, and $d$ is the greatest common divisor
of $r$ and $q$. Hence $a^{r'} - f_i b^{q'} = 0$ for some $f_i \in \overline{F}$. But then
$a^{r'}b^{-q'} \in F$ since it is a constant rational function defined over $F$
and $a^rb^{-q} = (a^{r'}b^{-q'})^d \in F$.   $\Box$\\

\textbf{Lemma on independence.} Let $\varphi$ be a homomorphism of
$A_1$ into $A_n$. Then the image of $A_1$ contains two elements with
algebraically independent leading forms.\\
\textbf{Proof.} Let $u = \varphi(x)$ and $v = \varphi(y)$ where $x$
and $y$ are generators of $A_1$ and let $B = F[u,v]$ be the image of $A_1$ in
$A_n$. If $\overline{u}$ and $\overline{v}$ are independent, we are done.
If not, then by Lemma on dependence
there exists a pair
of relatively prime natural numbers $(q, r)$ and $f \in F$ such that
$\overline{u^{q}} =
f\overline{v^{r}}$. Either $q$ or $r$ is not divisible by $p$. For arguments
sake assume that it is $q$.

We can find $k$ for which $kp + 1 \equiv 0 \pmod q$, $f_1 \in F$ and a
positive integer $s$ so that $\overline{u^{kp+1}} = f_1\overline{v^s}$.
Let
us replace the pair $(u, v)$ by the pair $(u_1 = u^{kp+1} - f_1v^s,
v_1 = v)$. The commutator $[u_1, v_1] = u^{kp}$ is a non-zero
element of the center $Z(B)$ of $B$.
If $\overline{u_1}$ and $\overline{v_1}$ are independent we are done,
otherwise repeat
the step above to get $(u_2, v_2)$, etc.. We claim that after a
finite number of steps we produce a pair of elements of $B$ with
independent leading forms.

Consider a function $\Def(a, b) = \deg(ab) -
\deg([a,b])$ on $A_n$. Let us check that
$\Def(u_{i+1}, v_{i+1}) < \Def(u_i,
v_i)$. We will do it for the first step since the computations are
the same for every step.

Since $\overline{u^{kp+1}} = f_1\overline{v^s}$ we see that
$\deg(u_1) < (kp + 1)\deg(u)$.
So $\Def(u_1, v_1) = \deg(u_1v_1) -
\deg([u^{kp+1} - f v^r, v]) = \deg(u_1) + \deg(v) - kp\deg(u) -
\deg([u, v])) < \deg(u) + \deg(v) - \deg([u, v] = \Def(u, v)$ since
$[u^{kp + 1} - f v^r, v] = u^{kp}[u, v]$ and $\deg(u_1) < (kp +
1)\deg(u)$.

By definition, $\Def(a, b) > 0$ if both $a$ and $b$ are not zero, so after at
most $\Def(u, v)$ steps we either get a pair with zero element or a
pair $U, \ V \in B$ with independent $U$ and $V$. Since $[u, v]
= 1$ the pair we start with does not contain zero. Similarly, since
$[u_i, v_i] \neq 0$ we cannot get a pair with
zero element. $\Box$\\

We can now see that $\gkdim(B) \geq 2$. Indeed, $U$ and $V$ are
``polynomials" of $u$ and $v$ and we may assume that the degrees of
these polynomials are at most $d$. Then the space of all polynomials
in $u, \ v$ of degree at most $N$ contains all polynomials in $U, \
V$ of degree at most $N\over d$. Since $\overline{U}$ and $\overline{V}$ are
algebraically independent the leading forms $\overline{U^i V^j}$ are
linearly independent over $F$ and hence $U^i V^j$ are linearly
independent over $F$. There are about $N^2\over 2d^2$ of these
monomials with $i + j \leq {N\over d}$ (exactly ${[{N\over d}] +
2\choose 2}$ where $[{N\over d}]$ is the integral part of $N\over
d$). So the dimension $d_N > {N^2\over 2d^2}$ and $\gkdim(B) \geq
2$.\\

\textbf{Theorem 2.} Let $\varphi$ be a homomorphism of $A_1$ into
$A_n$. Then $\varphi$ is an embedding. \\
\textbf{Proof.} Let $A_1 = F[x; y]$ and $u = \varphi(x), \ v =
\varphi(y)$. If $\varphi$ has a non-zero kernel take an element
$a$ in the kernel of smallest total degree possible. Since
both $[x, a]$ and $[y, a]$ are also in the kernel of $\varphi$ and
have smaller total degrees, $a \in Z(A_1) = F[x^p, y^p]$.
Therefore $u^p$ and $v^p$ are algebraically dependent and by Lemma
on GK-dimension $\gkdim(F[u^p, v^p]) \leq 1$ (recall that $u^p$
and $v^p$ commute). But $\gkdim(Z(B)) = \gkdim(B)$ for $B = F[u,
v]$. Indeed, $B = \sum u^iv^j Z(B)$ where $0 \leq i, \ j < p$ by
Lemma on center. So $d_N(B) \leq p^2 d_{N\over p}(Z(B))$ and
$d_N(B) \geq d_{N\over p}(Z(B))$. We showed above that $\gkdim(B)
\geq 2$. So $u^p$ and $v^p$ are algebraically independent and the
kernel of $\varphi$ consists of
zero only. $\Box$ \\

Theorem 2 cannot be extended to $A_2$.
Take $z = x_1 + y_1^{p-1}x_2, \  y_1, \  y_2$. Then $[z, y_1] = 1, \
[z, y_2] = y_1^{p-1},$ and $[z^p, y_2] =
\ad_{z}^p(y_2) = \ad_{z}^{p-1}(y_1^{p-1}) = (p-1)! = -1$. For $u_1 =
z + z^p y_1^{p-1}$, $v_1 = y_1$; $u_2 = y_2, \ v_2 = z^p$ the
commutation relations of $A_2$ are satisfied.
So $\phi$ which is defined by $\phi(x_i) = u_i$ and $\phi(y_i) = v_i$
is a homomorphism of $A_2$. If $B = \phi(A_2)$ then
$B = F[u_1, u_2; v_1, v_2] = F[z, y_1, y_2]$. Hence
$u_1^p \in Z(B)$ and $u_1^p \in
Z(F[z, y_1]) = F[z^p, y_1^p]$ since $u_1 \in F[z, y_1]$.
But $u_1^p$ should commute with
$u_2 = y_2$. Therefore $u_1^p \in F[z^{p^2}, y_1^p] = F[v_2^p, v_1^p]$,
$u_1^p, \ v_1^p, \ v_2^p$ are algebraically dependent, and $\phi$ has a non-zero kernel.

It is an exercise to check that $\gkdim(B) = 3$.
On the other hand, $\gkdim(A_2) = 4$ since $d_N = {N + 2n\choose 2n}$
for $A_n$ which gives $\gkdim(A_n) = 2n$.

A question about possible GK-dimensions of images of $A_n$ under
homomorphisms seems reasonable in this setting because clearly
the size of the kernel is large when the size of the image is small.
Say, $\gkdim(\varphi(A_n)) \leq 2n$ and if $\varphi$ is an injection
then $\gkdim(\varphi(A_n)) = 2n$. \\

\textbf{Theorem 3.} $\gkdim(\varphi(A_n))$ can be $n + i$ where $i = 1, 2,
\dots, n$ for a homomorphism $\varphi$ of $A_n$ into
$A_n$. \\
\textbf{Proof.} Denote $\varphi(A_n)$ by $B$ and by $u_i =
\varphi(x_i), \ v_i = \varphi(y_i)$. It is sufficient to show that
$\gkdim(B) = n + 1$ is possible for any $n$ because combining $u_1,
\ldots, u_m; v_1, \ldots, v_m$ of an appropriate map of $A_m$ to
$A_m$ with $x_{m+1}, \ldots, x_n; y_{m+1}, \ldots, y_n$ we will get
an image of $A_n$ of GK-dimension $m + 1 + 2(n - m)$.

Now we shall find $\varphi$ such that $\gkdim(B) = n + 1$.

Consider elements $z_{0,0} = 0, \ z_{m,0} = x_m - y_m^{p-1}z_{m-1,
0}$ for $m = 1, \ldots, n$, and $z_{m,i} = z_{m,0}^{p^i}$.
Then $[z_{i,0}, y_i] = 1$, $[z_{i-k,0}, y_i] = 0$  and $[z_{i-k,0},
x_i] = 0$ if $k > 0$, and $[z_{i,0}, z_{j,0}] = 0$. Therefore
$[z_{k,i}, z_{m,j}] = 0$.
We can get a relation between $z_{i, j}$ using the equality $(yx)^p - yx
= y^px^p$ if $[x, y] = 1$ (observe that $\ad_{yx} = \ad_{(yx)^p}$).
Take $y_mz_{m,0}  = y_mx_m - y_m^p z_{m-1,0}$. Then the summands in
the right side commute and $(y_mz_{m,0})^p = (y_m x_m)^p -
y_m^{p^2} z_{m-1,0}^p$. So $z_{m,0}^p = y_m^{-p}[(y_m x_m)^p -
y_m^{p^2}z_{m-1,0}^p - y_m x_m  + y_m^p z_{m-1,0}] =
y_m^{-p}[y_m^p x_m^p  - y_m^{p^2}z_{m-1,0}^p + y_m^p z_{m-1,0}] =
z_{m-1,0} - y_m^{p(p-1)}z_{m-1,0}^p + x_m^p$, i. e. $z_{m,1} =
z_{m-1,0} - y_m^{p(p-1)}z_{m-1,1} + x_m^p$. Since all summands in
the right side of this equality commute
$$z_{m,i+1} = z_{m-1,i} - y_m^{p^{i+1}(p-1)}z_{m-1,i+1} + x_m^{p^{i+1}}$$
and
$$z_{m-1,i} = z_{m,i+1} + y_m^{p^{i+1}(p-1)}z_{m-1,i+1} - x_m^{p^{i+1}}.$$

Then by induction we can prove that
$$z_{m-1,i} =  \sum_{k = 0}^{m-i-2}c_{m-1,i,k} z_{m, i+1+k} + c_{m-1,i}$$ and that
$$z_{m-j,i} =  \sum_{k = 0}^{m-i-j-1} c_{m-j,i,k}z_{m, i+j+k} + c_{m-j,i}$$
where $c_{i, j}, c_{i,j,k} \in Z(A_n)$. The sums are finite because
$z_{m, m} \in Z(A_n)$ (can be proved by
induction on $m$ starting with $z_{0,0} = 0$ since $z_{m,m} =
z_{m-1,m-1} - y_m^{p^{m}(p-1)}z_{m-1,m} + x_m^{p^{m}}$).

Now we can show that all $c_{i,j,k} \in F[y_1^p, \ldots, y_2^p]$.

Since $[z_{m,i},
y_j] = [z_{m-1,i-1}, y_j]  - y_m^{p^{i}(p-1)}[z_{m-1,i},  y_j]$,
we can deduce by induction on $m$ that
$[z_{m, i}, y_j]$ is zero if $j > m - i$, one if $j = m - i$, and
belongs to $F[y_1^p, \ldots, y_n^p]$ if $j < m - i$.
Since $z_{m-j,i} =  \sum_{k = 0}^{m-i-j-1} c_{m-j,i,k}z_{m, i+j+k} +
c_{m-j,i}$ we have $[z_{m-j,i}, y_l] = \sum_{k = 0}^{m-i-j-1}
c_{m-j,i,k}[z_{m, i+j+k}, y_l]$; this allows using $l = m - j - i,
\ldots, 1$ to check that all $c_{m-j,i,k} \in F[y_1^p, \ldots,
y_n^p]$.

All these computations were done to confirm that
$$z_{m,0} = \sum_{k=0}^{m-1} d_{m,k} z_{n, n-m+k} + d_m$$
where $d_m \in Z(A_n)$ and $d_{m,k} \in F[y_1^p, \ldots, y_n^p]$.

Recall that $x_m = z_{m,0} + y_m^{p-1}z_{m-1, 0}$ and so
$$x_m = \sum_{k=0}^{m-1} d_{m,k} z_{n, n-m+k} + d_m +
y_m^{p-1}(\sum_{k=0}^{m-2} d_{m-1,k} z_{n, n-m + 1+k} + d_{m-1}).$$
Therefore
$$u_m = x_m - d_m - y_m^{p-1}d_{m-1} \in B = F[y_1, \ldots, y_n; z_{n,0}].$$

Finally, $u_1, \ldots, u_n; y_1, \ldots, y_n$ define a homomorphism of
$A_n$ into $B$ and since any monomial in $B$ can be written as
$y_1^{j_1} \ldots  y_n^{j_n}z_{n,0}^i$, the Theorem is proved. $\Box$\\

By looking at $\varphi(F[x_1, \ldots, x_n])$ it is possible to show that
$\gkdim(\varphi(A_n)) \geq n$. This and Theorem 2 suggest the following\\

\textbf{Conjecture.} $\gkdim(\varphi(A_n)) > n$.\\

Acknowledgement. The author is grateful to the CRM (Centre de Recerca Matem\`{a}tica) in Bellaterra, Spain for the hospitality during the work on this project.

\begin{center}
{\bf References}
\end{center}

[AE] Adjamagbo, Pascal Kossivi; van den Essen, Arno A proof of the
equivalence of the Dixmier, Jacobian and Poisson conjectures. Acta
Math. Vietnam. 32 (2007), no. 2-3, 205--214.

[BCW] Bass, Hyman; Connell, Edwin H.; Wright, David The Jacobian
conjecture: reduction of degree and formal expansion of the inverse.
Bull. Amer. Math. Soc. (N.S.) 7 (1982), no. 2, 287--330.

[BK] Belov-Kanel, Alexei; Kontsevich, Maxim The Jacobian conjecture
is stably equivalent to the Dixmier conjecture. Mosc. Math. J. 7
(2007), no. 2, 209--218, 349.

[D] Dixmier, Jacques Sur les alg\`{e}bres de Weyl. (French) Bull. Soc.
Math. France 96 1968 209--242.

[GK] Gelfand, I. M.; Kirillov, A. A. Sur les corps li\'{e}s aux alg\`{e}bres
enveloppantes des alg\`{e}bres de Lie. (French) Inst. Hautes \'{E}tudes Sci.
Publ. Math. No. 31 1966 5--19.

[KL] Krause, G\"{u}nter R.; Lenagan, Thomas H. Growth of algebras and
Gelfand-Kirillov dimension. Revised edition. Graduate Studies in
Mathematics, 22. American Mathematical Society, Providence, RI,
2000. x+212 pp. ISBN: 0-8218-0859-1.

[T1] Tsuchimoto, Yoshifumi Endomorphisms of Weyl algebra and
p-curvatures, Osaka J. Math. 42 (2005), No. 2, 435–452.

[T2] Tsuchimoto, Yoshifumi Preliminaries on Dixmier conjecture. Mem.
Fac. Sci. Kochi Univ. Ser. A Math. 24 (2003), 43--59.

\end{document}